\documentclass[11pt]{article}

\usepackage{amsmath}
\usepackage{amssymb}
\usepackage{eucal}

\begin{document}

\baselineskip 16pt

\title{On finite $P\sigma T$-groups}

\author            
{ Alexander  N. Skiba \\
{\small Department of Mathematics,  Francisk Skorina Gomel State University,}\\
{\small Gomel 246019, Belarus}\\
{\small E-mail: alexander.skiba49@gmail.com}}

\date{}
\maketitle

\begin{abstract}      Let   $\sigma =\{\sigma_{i} | i\in I\}$ be some
 partition of the set of all primes $\Bbb{P}$ and  $G$ a  finite group. 
 $G$ is said to be \emph{$\sigma$-soluble} if every chief factor $H/K$ of 
$G$ is a
$\sigma _{i}$-group for some $i=i(H/K)$. 

A set
 ${\cal H}$ of subgroups of $G$  is said to be a  \emph{complete
Hall $\sigma $-set} of $G$  if every member $\ne 1$ of  ${\cal H}$
 is a Hall $\sigma _{i}$-subgroup of $G$ for some $\sigma _{i}\in \sigma $ and
   ${\cal H}$ contains exact one Hall  $\sigma _{i}$-subgroup of $G$ for every
 $i \in I$ such that $\sigma _{i}\cap   \pi (G)\ne \emptyset$.
A subgroup $A$ of $G$ is said to be   \emph{${\sigma}$-permutable} or
   \emph{${\sigma}$-quasinormal} in $G$
  if $G$ has a complete Hall $\sigma$-set $\cal H$  such that 
 $AH^{x}=H^{x}A$ for all $x\in G$ and all $H\in \cal H$.

We obtain  a characterization of  finite $\sigma$-soluble groups $G$ 
 in which  $\sigma$-quasinormality    is a transitive relation in $G$.

\end{abstract}

\footnotetext{Keywords: finite group, ${\sigma}$-quasinormal subgroup, $P\sigma T$-group, 
$\sigma$-soluble group, $\sigma$-nilpotent group.}

\footnotetext{Mathematics Subject Classification (2010): 20D10, 20D15, 20D30}
\let\thefootnote\thefootnoteorig

\section{Introduction}

Throughout this paper, all groups are finite and $G$ always denotes
a finite group. Moreover,  $\mathbb{P}$ is the set of all  primes,
  $\pi \subseteq  \Bbb{P}$ and  $\pi' =  \Bbb{P} \setminus \pi$. If
 $n$ is an integer, the symbol $\pi (n)$ denotes
 the set of all primes dividing $n$; as usual,  $\pi (G)=\pi (|G|)$, the set of all
  primes dividing the order of $G$.   
 $G$ is said to be  a $D_{\pi}$-group if $G$ possesses a Hall 
$\pi$-subgroup $E$ and every  $\pi$-subgroup of $G$ is contained in some 
conjugate of $E$.

In what follows, $\sigma$  is some partition of  
$\Bbb{P}$, that is,  $\sigma =\{\sigma_{i} |
 i\in I \}$, where   $\Bbb{P}=\cup_{i\in I} \sigma_{i}$
 and $\sigma_{i}\cap
\sigma_{j}= \emptyset  $ for all $i\ne j$;   $\Pi$ is always supposed
to be a    subset of the set $\sigma$ and $\Pi'= \sigma\setminus \Pi$.

By the analogy with the notation   $\pi (n)$, we write  $\sigma (n)$ to denote 
the set  $\{\sigma_{i} |\sigma_{i}\cap \pi (n)\ne 
 \emptyset  \}$;   $\sigma (G)=\sigma (|G|)$.  $G$ is said to be: 
\emph{$\sigma$-primary} \cite{1} if  $|\sigma (G)|\leq 1$; 
\emph{$\sigma$-decomposable} (Shemetkov \cite{Shem})  or 
\emph{$\sigma$-nilpotent} (Guo and Skiba  \cite{33}) if $G=G_{1}\times \dots \times G_{n}$ 
for some $\sigma$-primary groups $G_{1}, \ldots, G_{n}$; \emph{$\sigma$-soluble} \cite{1}
 if every chief factor of $G$ is $\sigma$-primary;   a \emph{$\sigma$-full group
 of Sylow type}  \cite{1} if every subgroup $E$ of $G$ is a $D_{\sigma _{i}}$-group for every
$\sigma _{i}\in \sigma (E)$.

A natural number $n$ is said to be a \emph{$\Pi$-number}
 if  $\sigma (n)\subseteq \Pi$.  A  subgroup $A$ of $G$ is said to be:
   a \emph{Hall $\Pi$-subgroup} of $G$ \cite{1, 2} if    $|A|$ is  a $\Pi$-number 
 and $|G:A|$ is a $\Pi'$-number; a \emph{$\sigma$-Hall subgroup} of $G$ if $A$
 is a Hall $\Pi$-subgroup of $G$ for some  $\Pi\subseteq \sigma$.

A     set  ${\cal H}$ of subgroups of $G$ is a
 \emph{complete Hall $\sigma $-set} of $G$ \cite{2, commun}  if
 every member $\ne 1$ of  ${\cal H}$ is a Hall $\sigma _{i}$-subgroup of $G$
 for some $\sigma _{i} \in \sigma$ and ${\cal H}$ contains exact one Hall
 $\sigma _{i}$-subgroup of $G$ for every  $\sigma _{i}\in  \sigma (G)$.

Recall  that a  subgroup $A$ of $G$ is said to be:    \emph{$\sigma$-permutable} or
 \emph{$\sigma$-quasinormal} in $G$ \cite{1} if $G$ possesses
 a complete Hall $\sigma$-set  ${\cal H}$ such that $AH^{x}=H^{x}A$ for all
  $H\in {\cal H}$ and  all $x\in G$; \emph{${\sigma}$-subnormal}
 in $G$ \cite{1}  if there is a subgroup chain  $$A=A_{0} \leq A_{1} \leq \cdots \leq
A_{t}=G$$  such that  either $A_{i-1}\trianglelefteq A_{i}$ or $A_{i}/(A_{i-1})_{A_{i}}$
 is  ${\sigma}$-primary for all $i=1, \ldots , t$.

In the classical case, when $\sigma =\sigma ^{0}=\{\{2\}, \{3\}, \ldots \}$, 
${\sigma}$-quasinormal subgroups are
 also called \emph{$S$-quasinormal} or  \emph{$S$-permutable} \cite{prod, GuoII}, and a 
 subgroup  $A$ of $G$ is subnormal in $G$ if and only if it is $\sigma 
^{0}$-subnormal in $G$. 

We say that  $G$ is a {\sl $P\sigma T$-group} \cite{1} if ${\sigma}$-quasinormality  
is a transitive relation in $G$, that is, if $K$ is a ${\sigma}$-quasinormal subgroup
 of $H$ and 
 $H$ is a ${\sigma}$-quasinormal subgroup of $G$, then  $K$ is a
 ${\sigma}$-quasinormal subgroup of $G$.   
 In the 
 case, when $\sigma =\{\{2\}, \{3\}, \ldots \}$, 
$P\sigma T$-groups are
called   \emph{$PST$-groups } \cite{prod}.

 In view of Theorem B in 
\cite{1}, $P\sigma T$-groups can be characterized as the groups in which every  
${\sigma}$-subnormal subgroup  is ${\sigma}$-quasinormal in $G$.

Our first  observation is the following fact, which generalizes the 
sufficiency condition in Theorem A of the paper \cite{1}. 

{\bf Theorem A.} {\sl   Let  $G$ have a normal $\sigma$-Hall  subgroup $D$ such that:

(i) {\sl $G/D$ is  a   $P\sigma T$-group, and  }

(ii) {\sl every  
$\sigma$-subnormal subgroup of $D$ is normal in $G$. }

{\sl  If $G$ is a $\sigma$-full group of Sylow type, then $G$ is a   
$P\sigma T$-group. }

}

{\bf Corollary 1.1} (See Theorem A in \cite{1}). {\sl  
 Let  $G$ have a normal $\sigma$-Hall  subgroup $D$ such that:   }

(i) {\sl $G/D$ is     $\sigma$-nilpotent, and  }

(ii) {\sl every  
 subgroup of $D$ is normal in $G$. }

{\sl  Then $G$ is a   
$P\sigma T$-group. }

In the case when $\sigma =\{\{2\}, \{3\}, \ldots \}$, we get from Theorem 
A the following 

{\bf Corollary 1.2} (See Theorem 2.4 in \cite{Agr}). {\sl 
  Let  $G$ have a normal Hall  subgroup $D$ such that:}

(i) {\sl $G/D$ is a   $PST$-group, and  }

(ii) {\sl every subnormal  
 subgroup of $D$ is normal in $G$. }

{\sl  Then $G$ is a   
$PST$-group. }

Recall that     $G^{\frak{N_{\sigma}}}$  denotes the \emph{$\sigma$-nilpotent
 residual} of $G$,
 that is,  the intersection of all normal subgroups $N$ of $G$ with 
$\sigma$-nilpotent quotient $G/N$;   $G^{\frak{N}}$
  denotes the \emph{nilpotent
 residual} of $G$ \cite{15}.

{\bf Definition 1.3.} We say that $G$ is a   
 \emph{special $P\sigma T$-group} provided  the $\sigma$-nilpotent 
residual $D=G^{{\frak{N}}_{\sigma}}$ of $G$ is contained in a Hall $\sigma 
_{i}$-subgroup $E$ of $G$ for some $i$  and   the  following  conditions hold: 

(i)  $D$ is a Hall subgroup of $G$ and every element of $G$ induces a 
power automorphism in $D$;

(ii) $D$ has a normal complement $S$ in $E$.

Note that if $G=C_{5}\times (C_{3}\rtimes C_{2})$, where $C_{3}\rtimes C_{2}\simeq S_{3}$
 and $\sigma =\{\{3, 5\}, \{3, 5\}'\}$, then $G$ is a special $P\sigma T$-group
with  $C_{3}=G^{\frak{N_{\sigma}}}$.

The following   theorem shows that every special 
$P\sigma T$-group is a  $P\sigma T$-group.

{\bf Theorem B.} {\sl   Suppose that $G$ has a $\sigma$-nilpotent normal
Hall  subgroup $D$ with $\sigma$-nilpotent 
quotient $G/D$ such that $G/ O^{\sigma _{i}}(D)$ is a special $P\sigma T$-group 
 for each $\sigma _{i} \in \sigma (D) $. Then $G$ is a   
$P\sigma T$-group. }

Generalizing the concept of  complete Wielandt $\sigma$-set 
 of a group in  \cite{33}, we say that 
 a complete Hall $\sigma $-set   $\cal H$ of $G$
is a \emph{generalized Wielandt $\sigma$-set} of $G$  if  every member $H$ of
 $\cal H$  is $\pi (G^{{\frak{N}}_{\sigma}})$-supersoluble.

 Using  Theorem B, we prove also the  following revised version of Theorem A in 
\cite{1}.

{\bf Theorem C. } {\sl Let $G$ be    $\sigma$-soluble
 and $D=G^{\frak{N_{\sigma}}}$. Suppose that 
 $G$ has a generalized Wielandt $\sigma$-set. 
   Then  $G$ is a $P\sigma 
T$-group if and only if  the following conditions hold:}   

(i) {\sl $D$ is an abelian  Hall subgroup of $G$ of odd
 order and  every element of $G$ induces a
 power automorphism in $D$;  }

(ii) {\sl $G/ O^{\sigma _{i}}(D)$ is a special $P\sigma T$-group 
 for each $\sigma _{i} \in \sigma (D) $.}

{\bf Corollary 1.4} (See Theorem 2.3 in  \cite{Agr}). {\sl Let 
$G$ be a soluble and $D=G^{\frak{N}}$.  If  $G$ is a $PST$-group, then  $D$  is  an
abelian  Hall subgroup of $G$ of odd  order and every element
 of $G$ induces a power automorphism 
in  $D$.  }

\section{Some preliminary  results}

In view of Theorems A and B in \cite{2}, the following fact is true.  

{\bf Lemma 2.1.} {\sl If $G$ is $\sigma$-soluble, then $G$ is a $\sigma$-full group
 of Sylow type.    
}

We use ${\mathfrak{N}}_{\sigma}$
  to denote the class  of all  $\sigma$-nilpotent groups.

{\bf Lemma 2.2 } (See  Corollary 2.4 and Lemma 2.5  in \cite{1}).  {\sl  The class 
 ${\mathfrak{N}}_{\sigma}$               is closed under taking  
products of normal subgroups, homomorphic images and  subgroups. Moreover, if  $E$ is a normal 
subgroup of $G$ and  $E/E\cap \Phi (G)$ is $\sigma$-nilpotent, then 
$E$ is $\sigma$-nilpotent.    }

In view of Proposition 2.2.8  in \cite{15},  
 we get from Lemma 2.2 the following

{\bf Lemma 2.3.}    {\sl If 
$N$ is a normal subgroup of $G$, then
 $$(G/N)^{{\frak{N}}_{\sigma}}=G^{{\frak{N}}_{\sigma}}N/N.$$  }

{\bf Lemma 2.4} (See  Knyagina and  Monakhov \cite{knyag}). {\sl
Let $H$, $K$  and $N$ be pairwise permutable
subgroups of $G$ and  $H$ be  a Hall subgroup of $G$. Then $$N\cap HK=(N\cap H)(N\cap K).$$}

{\bf Lemma 2.5.}   {\sl  The following statements hold:}

(i) {\sl $G$ is a 
 $P\sigma T$-group if and only if every  $\sigma$-subnormal subgroup of 
$G$ is $\sigma$-quasinormal in $G$. }

(ii) {\sl 
If  $G$ is a 
 $P\sigma T$-group, then every   quotient $G/N$ of $G$ is also a   
$P\sigma T$-group. }

(iii) {\sl     If  $G$ is a  special 
 $P\sigma T$-group, then every   quotient $G/N$ of $G$ is also a  special  
$P\sigma T$-group. }

 {\bf Proof.} (i) This follows from the fact (see  Theorem  B in \cite{1}) that
 every $\sigma$-quasinormal 
subgroup of $G$
 is $\sigma$-subnormal in $G$. 

(ii) Let $H/N$ be a $\sigma$-subnormal subgroup of $G/N$. Then 
$H$   is a $\sigma$-subnormal subgroup of $G$ by Lemma 2.6(5) in \cite{1}, so $H$ is
 $\sigma$-quasinormal in $G$ by hypothesis and Part (i). Hence $H/N$  is
  $\sigma$-quasinormal in 
$G/N$ by Lemma 2.8(2) in \cite{1}. Hence   $G/N$  is  a   
$P\sigma T$-group by Part (i).

(iii)  Suppose that $D=G^{{\frak{N}}_{\sigma}}$ is a Hall subgroup of $G$ and 
 $D\leq E$, where $E=D\times S$ is a Hall $\sigma 
_{i}$-subgroup $E$ of $G$,  and every element of $G$ induces a power 
automorphism in $D$. 
 Then     $EN/N$ is a  Hall $\sigma 
_{i}$-subgroup  of $G/N$   and  
 $DN/N=(G/N)^{{\frak{N}}_{\sigma}}$ is a Hall subgroup of $G/N$  by Lemma 2.3. 
Moreover,  $EN/N=(DN/N)(SN/N)$ and,  by Lemma 2.4, 
 $$DN\cap SN=N(D\cap SN)= N(D\cap S)(D\cap N)=N(D\cap N)=N,$$ which implies that 
$(DN/N)\cap (SN/N)=1$. Hence $EN/N=(DN/N)\times (SN/N)$.

 Finally, let $H/N\leq DN/N$. Then
 $H=N(H\cap D)$, where $H\cap D$ is normal in $G$ by hypothesis. But then 
$H/N= N(H\cap D)/N$ is normal in $G/N$, so every element of $G/N$ induces a power 
automorphism on $DN/D$.  Hence $G/N$ is a  special  
$P\sigma T$-group.

The lemma is proved.

\section{Proofs of the results}

{\bf Proof of Theorem A.}  Since  $G$ is a $\sigma$-full group of Sylow type by hypothesis,
 it  possesses a 
complete Hall $\sigma $-set   ${\cal H}=\{H_{1}, \ldots , H_{t} \}$, and a
 subgroup $H$ of $G$ is 
$\sigma$-quasinormal in $G$ if and only if $HH_{i}^{x}=H_{i}^{x}H$ for all
 $H_{i} \in {\cal H}$ and $x\in G$.  
  We can assume without loss  of generality  that $H_{i}$ is a 
$\sigma _{i}$-group for all $i=1, \ldots , t$.

  Assume that
this theorem  is false and let $G$ be a counterexample of minimal order.  
Then  $D\ne 1$ and   for some  $\sigma$-subnormal subgroup $H$ of $G$ and for
some $x\in G$ and $k\in I$
 we have  $HH_{k}^{x}\ne H_{k}^{x} H$ by Lemma 2.5(i).  Let $E=H_{k}^{x}$.

(1) {\sl The hypothesis holds for every quotient  $G/N$ of $G$. }

 It
 is clear that  $G/N$  is a $\sigma$-full group of Sylow type and 
  $DN/N$   is  
a normal $\sigma$-Hall  subgroup of $G/N$.  On the other hand,  
$$(G/N)/(DN/N)  \simeq G/DN \simeq (G/D)/(DN/D),$$ so $(G/N)/(DN/N)$ is a  
$P\sigma T$-group  by Lemma 2.5(ii).  Finally, let $H/N$ be a  $\sigma$-subnormal 
subgroup of $DN/N$.     
Then  $H=N(H\cap D)$  and, by Lemma 2.6(5) in \cite{1},  $H$ is $\sigma$-subnormal in $G$. 
  Hence $H\cap D$ is 
$\sigma$-subnormal in $D$ by Lemma 2.6(1) in \cite{1}, so $H\cap D$ is 
normal in $G$ by   hypothesis. Thus $H/N=N(H\cap D)/N$  is normal in 
$G/N$.   Therefore the  hypothesis holds on  $G/N$.

(2) {\sl $H_{G}=1$. }

Assume that $H_{G}\ne 1$. Clearly,   $H/H_{G}$ is  $\sigma$-subnormal in 
$G/H_{G}$. 
Claim (1) implies that  the   hypothesis holds  for 
$G/H_{G}$, so the choice of $G$ implies that  $G/H_{G}$ is a $P\sigma T$-group. Hence 
$$(H/H_{G})(EH_{G}/H_{G})=
 (EH_{G}/H_{G})(H/H_{G}).$$ by Lemma   2.5(i). Therefore $EH=EHH_{G}$ is a subgroup of $G$
 and so $HE=EH$,  a contradiction. 
 Hence $H_{G}= 1$.

(3) {\sl $DH=D\times H$. }

By Lemma 2.6(1) in \cite{1}, $H\cap D$ is $\sigma$-subnormal in  $D$. Hence  
$H\cap D$ is normal in $G$ by hypothesis, which implies that $H\cap D=1$ by Claim  (2).  
Lemma 2.6(1) in \cite{1} implies also that   $H$ is $\sigma$-subnormal in $DH$. But
 $H$ is a $\sigma$-Hall subgroup of $DH$  since $D$ is a 
$\sigma$-Hall subgroup  of $G$ and $H\cap D=1$.   Therefore $H$ is normal in $DH$ by Lemma 
2.6(10) in \cite{1}, so $DH=D\times H$.

 {\sl Final contradiction.} Since $D$ is a $\sigma$-Hall subgroup of $G$, then either $E\leq D$
 or $E\cap D=1$. But the former
 case is impossible by Claim (3) since $HE\ne EH$, so  $E\cap D=1$.  Therefore $E$ is a 
$\Pi'$-subgroup of $G$, where $\Pi =\sigma (D)$.  
 By the 
Schur-Zassenhaus theorem, $D$ has a complement $M$ in $G$.  Then $M$ is a 
Hall  $\Pi'$-subgroup of $G $ and so for some $x\in G$ we have $E\leq 
M^{x}$ since $G$ is a $\sigma$-full group of Sylow type.  On the other hand, $H\cap M^{x}$
 is a Hall  $\Pi'$-subgroup of 
$H$ by Lemma 2.6(7) in \cite{1} and hence $H\cap M^{x}=H\leq  M^{x}$ since $H\cap D=1$ by
 Claim 
(3).  Lemma 2.6(1) 
in \cite{1} implies that $H$ is  $\sigma$-subnormal in $M^{x}$.  But $M^{x}\simeq 
G/D$ is a $P\sigma T$-group by hypothesis, so $HE=EH$ by Lemma 2.5(i). 
This contradiction completes the proof of the theorem.

{\bf Lemma 3.1.}   {\sl   If    $G$ is a  special 
 $P\sigma T$-group, then it is   a $P\sigma T$-group. }

 {\bf Proof.}  Let $D=G^{{\frak{N}}_{\sigma}}$ and $E$ be a normal Hall $\sigma 
_{i}$-subgroup 
 of $G$ such that $E=D\times S$.   Since $G/D$ is 
 $\sigma$-nilpotent, $G$ is $\sigma$-soluble. Hence 
   $G$ is a 
$\sigma$-full group of Sylow type by Lemma  2.1. Therefore $G$ possesses a 
complete Hall $\sigma $-set   ${\cal H}=\{H_{1}, \ldots , H_{t} \}$, and a
 subgroup $H$ of $G$ is 
$\sigma$-quasinormal in $G$ if and only if $HH_{j}^{x}=H_{j}^{x}H$ for all
 $H_{j} \in {\cal H}$ and $x\in G$.
  We can assume without loss  of generality  that $H_{j}$ is a 
$\sigma _{j}$-group for all $j=1, \ldots , t$.

  Assume that this lemma is false and let $G$ be a counterexample of minimal 
order. Then  $G$ is not $\sigma$-nilpotent, and 
 for some  $\sigma$-subnormal subgroup $H$ of $G$ and for
some $x\in G$ and $k\in I$
 we have  $HH_{k}^{x}\ne H_{k}^{x} H$ by Lemma 2.5(i).  Let $E=H_{k}^{x}$.  The subgroup 
 $S$ is
 normal in $G$ since it
 is characteristic in $E$. Since $G$ is not $\sigma$-nilpotent, $D\ne 1$. On the other hand,
 Theorem A  and the choice of $G$ 
imply that $S\ne 1$ since every subgroup of $D$ is normal in $G$ by hypothesis.
 Let $R$ and $N$ be minimal normal subgroups of $G$ such that  $R\leq 
D$ and $N\leq S$.  Then $R$ is a group of order $p$ for some prime $p$.
 Hence $R\cap HN\leq O_{p}(HN)\leq P$, where $P$ is a Sylow 
$p$-subgroup of $H$ since $\pi (D)\cap \pi (S)=\emptyset$, so $R\cap HN=R\cap H$. 

    The hypothesis holds for $G/R$ and 
$G/N$ by Lemma 2.5(iii). Hence the choice of $G$ and  Lemma 2.5(i) imply that
$$EHR/R=(ER/R)(HR/R)=(HR/R)(EHR/R)$$ 
 and  so $EHR$ is a subgroup of $G$.  
 Similarly we get that  $EHN$ is a 
subgroup of $G$. Since $|R|=p$ and $EH$ is not a subgroup of $G$, $R\cap 
E=1$.  Therefore from Lemma 2.4 we get that  that $R\cap EHN=R\cap E(HN)=(R\cap E)(R\cap HN)=R\cap HN$. 
 Hence  $$EHR\cap EHN=E(HR\cap EHN)=EH(R\cap EHN)=
EH(R\cap HN)=$$$$=EH(R\cap HN)=EH(R\cap H)=EH$$  is a subgroup of $G$. Hence  
$HE=EH$, a contradiction. The lemma  is proved.

{\bf Lemma 3.2. } {\sl  If ${\cal H}=\{H_{1}, \ldots , H_{t} \}$ is   
  a  generalized Wielandt $\sigma$-set of $G$, then 
 $${\cal H}_{0}=\{H_{1}N/N, \ldots , H_{t}N/N \}$$ is   
  a  generalized Wielandt $\sigma$-set of $G/N$.}

{\bf Proof.} It is clear that ${\cal H}_{0}$ is a complete Hall $\sigma$-set of $G/N$. Now let $D=G^{{\frak{N}}_{\sigma}}$ and $\pi =\pi (G^{{\frak{N}}_{\sigma}})$. 
Then $(G/N)^{{\frak{N}}_{\sigma}}=DN/N$ by Lemma 2.3, so
 $$\pi _{0}=\pi ((G/N)^{{\frak{N}}_{\sigma}})
=\pi (DN/N)\subseteq  \pi (D)=\pi.$$  Hence every member $H_{i}$ of ${\cal H}$ is  $\pi 
_{0}$-supersoluble, so $H_{i}N/N$ is $\pi _{0}$-supersoluble.  Hence  ${\cal H}_{0}$ is   
  a  generalized Wielandt $\sigma$-set of $G/N$.          
The lemma is proved. 

{\bf Proof of Theorem B.}  Clearly, $G$ is  $\sigma$-soluble, so $G$ is a 
$\sigma$-full group of Sylow type by Lemma  2.1. Therefore $G$ possesses a 
complete Hall $\sigma $-set   ${\cal H}=\{H_{1}, \ldots , H_{t} \}$, and a
 subgroup $H$ of $G$ is 
$\sigma$-quasinormal in $G$ if and only if $HH_{i}^{x}=H_{i}^{x}H$ for all
 $H_{i} \in {\cal H}$ and $x\in G$.
  We can assume without loss  of generality  that $H_{i}$ is a 
$\sigma _{i}$-group for all $i=1, \ldots , t$.

  Assume that
this theorem  is false and let $G$ be a counterexample of minimal order.
 Then  $D\ne 1$ and for some  $\sigma$-subnormal subgroup $H$ of $G$ and for
some $x\in G$ and $k\in I$
 we have  $HH_{k}^{x}\ne H_{k}^{x} H$  by Lemma 2.5(i). Let $E=H_{k}^{x}$.

(1) {\sl $G$ is not  a special $P\sigma T$-group} (This follows from Lemma 
3.1 and the choice of $G$).

(2) $|\sigma (D)| > 1$.

Indeed, suppose that $\sigma (D)=\{\sigma _{i}\}$.  Then  $O^{\sigma 
_{i}}(D)=1$, so $G\simeq  G/O^{\sigma 
_{i}}(D)$ is a special $P\sigma T$-group by hypothesis, contrary to Claim (1).

(3) {\sl The hypothesis holds for every quotient  $G/N$ of $G$, where $N\leq D$.}

First we show that $(G/N)/ O^{\sigma _{i}}(DN/N)$ is a special $P\sigma T$-group 
 for each $\sigma _{i} \in \sigma (DN/N) $. Note that $\sigma _{i} 
\in \sigma (DN/N)=\sigma (D/(D\cap N))\subseteq \sigma (D)$, so $G/ O^{\sigma _{i}}(D)$ is a special
 $P\sigma T$-group  by hypothesis.   It is not difficult to show that  $$O^{\sigma 
_{i}}(D)N/N=O^{\sigma _{i}}(D/N).$$ Hence $$(G/N)/(O^{\sigma _{i}}(D/N))=(G/N)/(O^{\sigma 
_{i}}(D)N/N)\simeq G/NO^{\sigma 
_{i}}(D)\simeq$$ $$\simeq  (G/O^{\sigma 
_{i}}(D))/(O^{\sigma 
_{i}}(D)N/O^{\sigma _{i}}(D))$$ is a special   $P\sigma T$-group  by Lemma 
2.5(iii).

 It is clear also that $DN/N\simeq D/D\cap N$ is a $\sigma$-nilpotent
 normal Hall subgroup of $G/N$ 
with $\sigma$-nilpotent quotient $$ (G/N)/(DN/N)\simeq G/DN\simeq 
(G/D)/(DN/D)$$ by Lemma 2.2. Hence we have (3).

(4) {\sl If $N$ is a minimal normal subgroup of $G$ contained in $D$, then $EHN$ is
 a subgroup of $G$.}

Claim (3) and the choice of $G$ implies that  the conclusion of the theorem holds  for
 $G/N$. On the other hand, $EN/E$ is a Hall $\sigma _{k}$-subgroup of $G/N$ and, by Lemma 2.6(4) in \cite{1}, 
    $HN/N$ is a $\sigma$-subnormal subgroup of 
$G$.   Note also that $G/N$ is $\sigma$-soluble, so every two Hall 
$\sigma _{k}$-subgroups of $G/N$ are conjugate by Lemma 2.1. Thus,
 $$(HN/N)(EN/N)=(EN/N)(HN/N)=EHN/N$$
 by Lemma 2.5(i).   Hence $EHN$ is
 a subgroup of $G$.

{\sl Final contradiction.} Since   $|\sigma (D)| > 1$ by Claim (2) and $D$
 is $\sigma$-nilpotent, $G$  
has at least two  $\sigma$-primary  minimal normal subgroups $R$ and $N$ such that $R, N\leq D$ 
and $\sigma (R)\ne \sigma (N)$.  Then  at least one of the subgroups 
$R$ or $N$, $R$ say, is a $\sigma _{i}$-group for some $i\ne k$. Then  
 $R\cap HN\leq O_{\sigma _{i}}(HN)\leq V$, where $V$ is a Hall $\sigma _{i}$-subgroup of $H$,
 since $N$ is a $\sigma _{i}'$-group and $G$ is a $\sigma$-full group of 
Sylow type.  Hence  $R\cap HN=R\cap H$.
Claim (4) implies that   $EHR$ and $EHN$ are subgroups of $G$. 
Now,  arguing similarly as in the proof of Lemma 3.1, one can show that 
  $EHR\cap EHN=EH$ is
 a subgroup of $G$, so $HE=EH$. 
This contradiction  
 completes the proof of the result.

{\bf Proof of Theorem C.}  
 Let $\pi =\pi (D)$ and  ${\cal H}=\{H_{1}, \ldots , H_{t} \}$ be  
  a  generalized Wielandt $\sigma$-set of $G$.
  We can assume without loss  of generality  that $H_{i}$ is a 
$\sigma _{i}$-group for all $i=1, \ldots , t$.   Since $G$ is 
$\sigma$-soluble by hypothesis, $G$ is a $\sigma$-full group of Sylow type by Lemma 2.1.

 {\sl Necessity.}   Assume that
this    is false and let $G$ be a counterexample of minimal order.  
Then $D\ne 1.$

(1) {\sl If $R$ is a non-identity normal subgroup of $G$, then the 
 hypothesis holds for $G/R$. Hence the necessity condition of the theorem  holds for $G/R$}
(Since the hypothesis holds for $G/R$ by Lemmas 2.5(ii) and 3.2, this follows from 
 the choice of $G$).

(2) {\sl If $E$ is a proper $\sigma$-subnormal subgroup of $G$, then 
 $E^{\frak{N_{\sigma}}}\leq D$ and the necessity condition of the theorem holds for $E$.   }

Every  $\sigma$-subnormal   subgroup $H$ of $E$ is $\sigma$-subnormal 
in $G$ by Lemma 2.6(2) in  \cite{1} and hence   $H$ is  $\sigma$-quasinormal 
in $G$ by hypothesis and Lemma 2.5(i). Thus $H$ is $\sigma$-quasinormal in $E$ by Lemma 
2.8(1)  in  \cite{1} since $G$ is a $\sigma$-full group 
  of Sylow type.  Thus, $E$ is a $\sigma$-soluble  
$P\sigma T$-group.   
It is clear that
 $E$ possesses a complete Hall $\sigma$-set
 ${H}_{0}= \{E_{1}, \ldots , E_{n} \}$ such  that $E_{i}\leq H_{i}^{x_{i}}$ for 
some $x_{i}\in G$ for all $i=1, \ldots , n.$   Hence every member of 
${H}_{0}$ is $\pi$-supersoluble.
Moreover,  since $$E/E\cap D\simeq ED/D\in {\frak{N_{\sigma}}}$$ 
    and ${\frak{N_{\sigma}}}$ is a hereditary  class by Lemma 2.2, we have   
$E/E\cap D    \in {\frak{N_{\sigma}}}$. Hence $E^{\frak{N_{\sigma}}}\leq 
E\cap D$.  Therefore, $\pi _{0}=\pi (E^{\frak{N_{\sigma}}})\subseteq \pi$. 
Hence  every member of ${H}_{0}$ is $\pi _{0}$-supersoluble. Hence 
${H}_{0}$ is a    generalized Wielandt $\sigma$-set of $E$.

Therefore the hypothesis holds for $E$, 
so the necessity condition of the theorem holds for $E$ by the choice of $G$.

(3) {\sl $D$ is nilpotent.}

 Assume that this is false and let $R$ be
  a minimal   normal subgroup of $G$. Then
$RD/R=(G/R)^{{\frak{N_{\sigma}}}}$ is abelian by  Lemma 2.3 and Claim (1).
 Therefore $R\leq D$, $R$   
is the unique minimal   normal subgroup of $G$ and $R\nleq \Phi (G)$ 
by  Lemma 2.2.   Let $V$ be a maximal 
subgroup of $R$. 
Since $G$ is $\sigma$-soluble by hypothesis, $R$  
 is a $\sigma _{i}$-group for some $i$.   Hence $V$ is $\sigma$-subnormal in $G$ by Lemma 
2.6(6)  in \cite{1}, so $V$ is  $\sigma$-quasinormal in $G$  by hypothesis and Lemma 2.5(i).
 Then  
$R\leq D\leq  O^{\sigma _{i}}(G)\leq N_{G}(V)$  by Lemma 3.1 in \cite{1}.
  Hence $R$ is abelian, so   $R=C_{G}(R)$ is a $p$-group for some prime 
$p$ by \cite[A, 15.2]{DH}.    

 It is clear that $R\leq H_{i}\cap D$ for some $i$.
 Then  $H_{i}$ is $p$-supersoluble by hypothesis, so some subgroup $L$ of $R$ of order $p$ 
is normal in $H_{i}$.  On the other hand,  $L$ is clearly  $\sigma$-quasinormal in
 $G$ and hence  $G=H_{i}O^{\sigma 
_{i}}(G)\leq N_{G}(L)$ by Lemma 3.1 in \cite{1}, so $R=L$.   
Therefore  $G/C_{G}(R)=G/R$ is a cyclic group. Hence
 $G$ is supersoluble and therefore  $D$ is nilpotent.

(4) {\sl  $D$ is a Hall subgroup of $G$. }

 Suppose
that this is false and let $P$ be a  Sylow $p$-subgroup of $D$ such
that $1 < P < G_{p}$, where $G_{p}\in \text{Syl}_{p}(G)$.  We can assume 
without loss of generality that $G_{p}\leq H_{1}$.

(a)  {\sl    $D=P$ is  a minimal normal subgroup of $G$. }

Let $R$ be a minimal normal subgroup of $G$ contained in $D$. 
 Since
 $D$ is  nilpotent by Claim (3),   $R$ is a $q$-group    for some prime   
$q$. Moreover, 
$D/R=(G/R)^{\mathfrak{N}_{\sigma}}$  is a Hall subgroup of $G/R$ by
Claim (1) and Lemma 2.3.  Suppose that  $PR/R \ne 1$. Then  $PR/R \in \text{Syl}_{p}(G/R)$. 
If $q\ne p$, then    $P \in \text{Syl}_{p}(G)$. This contradicts the fact 
that $P < G_{p}$.  Hence $q=p$, so $R\leq P$ and therefore $P/R \in 
\text{Syl}_{p}(G/R)$ and we again get  that 
$P \in \text{Syl}_{p}(G)$. This contradiction shows that  $PR/R=1$, which implies that 
  $R=P$  is the unique minimal normal subgroup of $G$ contained in $D$. Since $D$ is nilpotent,
 a $p'$-complement $E$ of $D$ is characteristic in 
$D$ and so it is normal in $G$. Hence $E=1$, which implies that $R=D=P$.

(b) {\sl $D\nleq \Phi (G)$.    Hence for some maximal subgroup
 $M$ of $G$ we have $G=D\rtimes M$  }  (This follows  from  Lemma 2.2 since $G$
 is not $\sigma$-nilpotent).

(c) {\sl If $G$ has a minimal normal subgroup $L\ne D$, then $G_{p}=D\times (L\cap G_{p})$.
  Hence $O_{p'}(G)=1$. }

Indeed, $DL/L\simeq D$ is a Hall 
subgroup of $G/L$ by Claim (1). Hence  $G_{p}L/L=RL/L$, so $G_{p}=D\times (L\cap G_{p})$.
 Thus  $O_{p'}(G)=1$ since $D < G_{p}$ by Claim (a).

(d)  {\sl   $V=C_{G}(D)\cap M$ is a  normal subgroup of $G$ and 
 $C_{G}(D)=D\times V \leq H_{1}$.  }

In view of  Claim  (b),  $C_{G}(D)=D\times V$, where $V=C_{G}(D)\cap M$ 
is a normal  subgroup of $G$. By Claim (a), $V\cap D=1$ and 
hence $V\simeq DV/D$ is $\sigma $-nilpotent by Lemma 2.2.  Let $W$ be a $\sigma 
_{1}$-complement of $V$. Then $W$  is characteristic in $V$ and so it is normal 
in $G$.    Therefore we have  (d) by Claim (c).

(e)  $G_{p}\ne H_{1}$.

Assume that $G_{p}=H_{1}$.  Let $Z$ be a subgroup of order $p$ in $Z(G_{p})\cap D$.
Then, since $ D\leq O^{\sigma _{1}}(G)=O^{p}(G)$, $Z$ is normal in 
$G$ by Lemma 3.1 in \cite{1}. Hence $D=Z < G_{p}$ and so $D < C_{G}(D)$. 
   Then  $V=C_{G}(D)\cap M\ne 1$ is a normal subgroup of $G$ and 
  $V\leq H_{1}=G_{p}$ by Claim (d). Let $L$ be a   minimal
 normal subgroup of $G$ contained in $V$. Then  $G_{p}=D\times L$ is a normal  
elementary abelian subgroup of $G$. 
  Therefore every subgroup of $G_{p}$ is 
normal in $G$ by Lemma 3.1  in \cite{1}.  Hence  $|D|=|L|=p$.
Let $D=\langle a \rangle$,  $L=\langle b \rangle$ and $N=\langle ab \rangle$.  
Then $N\nleq D$, so in view of the $G$-isomorphisms
 $$DN/D\simeq N\simeq NL/L= G_{p}/L=DL/L\simeq D $$  we get that 
$G/C_{G}(D)=G/C_{G}(N)$ is a $p$-group since $G/D$ is $\sigma$-nilpotent by Lemma 
2.2.
But then Claim (d) implies that  $G$ is a $p$-group. This 
contradiction shows that we have (e).

{\sl Final contradiction for (4).} In view of Theorem A in \cite{2}, $G$ has a $\sigma 
_{1}$-complement $E$ such that $EG_{p}=G_{p}E$. 
Let $V=(EG_{p})^{{\frak{N}}_{\sigma}}$.  By 
Claim (e), $EG_{p}\ne G$.     On the other hand, since $   D\leq 
EG_{p}$ by Claim (a),  $EG_{p}$ is $\sigma$-subnormal in $G$ by Lemma 
2.6(5)  in  \cite{1}. 
 Therefore   the necessity condition of the theorem holds for $EG_{p}$  by Claim (2). 
Hence 
   $V$  is a Hall subgroup of $EG_{p}$.  Moreover,  by Claim (2) we have 
$V\leq D$,  so  for a Sylow 
$p$-subgroup $V_{p}$ of $V$ we have $|V_{p}|\leq |P| < |G_{p}|$. 
Hence  $V$ is a $p'$-group and so   
 $V\leq C_{G}(D)\leq H_{1}=G_{p}$. Thus   $V=1$.  
Therefore $EG_{p}=E\times G_{p}$ is $\sigma$-nilpotent and so $E\leq C_{G}(D)\leq     
H_{1}$ by Claim (d). Hence $E=1$ and  so $ D =1$, a contradiction.  Thus,   
$D$ is a Hall subgroup of $G$. 

(5)  {\sl $G/ O^{\sigma _{i}}(D)$ is a special $P\sigma T$-group 
 for each $\sigma _{i} \in \sigma (D) $.}

First assume that $O^{\sigma 
_{i}}(D)\ne 1$ and let $N$ be a minimal normal subgroup of $G$ contained 
in $ O^{\sigma _{i}}(D)$. Then $G/N$ is a $P\sigma T$-group by  Lemma 2.5(ii), so     
the choice of $G$ implies that   
  $$(G/N)/ O^{\sigma _{i}}(D/N)=
(G/N)/ (O^{\sigma _{i}}(D)/N)\simeq G/ O^{\sigma _{i}}(D)$$  is a special 
$P\sigma T$-group.
Now assume  that  $O^{\sigma _{i}}(D)=1$, 
that is, $D$ is a $\sigma _{i}$-group. Since $G/D$ is  $\sigma$-nilpotent 
by Lemma 2.2, $H_{i}/D$ is normal in $G/D$ and hence $H_{i}$ is normal in 
$G$.   Therefore all subgroups of $D$ are 
$\sigma$-permutable  in $G$ by Lemma 2.3(2)(3) and hypothesis. 
Since $D$ is a normal Hall subgroup of $H_{i}$, it has a 
complement $S$ in $H_{i}$ by the Schur-Zassenhaus theorem. Lemma 3.1 in 
\cite{1} implies that $D\leq   O^{\sigma _{i}}(G)\leq N_{G}(S)$. Hence 
$H_{i}=D\times S$. Therefore $$G=H_{i}O^{\sigma _{i}}(G)=SO^{\sigma _{i}}(G)\leq N_{G}(L)$$ for 
every subgroup $L$ of $D$. Hence every element of $G$ induces a power 
automorphism in $D$. Hence $G$ is a special $P\sigma T$-group.

(6) {\sl Every subgroup $H$ of $D$ is normal in $G$. Hence every element of
 $G$ induces a power automorphism in $D$. }

Since $D$ is  nilpotent by Claim (3), it is enough
 to consider 
the case when $H$ is a subgroup of the Sylow $p$-subgroup $P$ of $D$ for  some prime $p$.
For some $i$ we have   $P\leq O_{\sigma _{i}}(D)=H_{i}\cap D$. On the other hand, we have
 $$D=O_{\sigma _{i}}(D)\times O^{\sigma _{i}}(D)$$ since 
$D$ is nilpotent.  Then  $$HO^{\sigma _{i}}(D)/O^{\sigma _{i}}(D)\leq
 D/O^{\sigma _{i}}(D)=(G/O^{\sigma _{i}}(D))^{{\frak{N}}_{\sigma}},$$ so 
 $HO^{\sigma _{i}}(D)/O^{\sigma _{i}}(D)$  is 
normal in $G/O^{\sigma _{i}}(D)$ by Claim (5). Hence  $HO^{\sigma 
_{i}}(D)$ is normal in $G$, which implies that $$H=H(O^{\sigma _{i}}(D)\cap 
O_{\sigma _{i}}(D))=HO^{\sigma _{i}}(D)\cap 
O_{\sigma _{i}}(D)$$ is normal in $G$.

 (7) {\sl  If  $p$ is a  prime such that $(p-1, |G|)=1$, then  $p$
does not divide $|D|$. In particular,  $|D|$ is odd. }

Assume that this is false.
 Then, by Claim (6),  $D$ has a maximal subgroup $E$ such that
$|D:E|=p$ and  $E$ is normal in $G$. It follows that  $C_{G}(D/E)=G$ since $(p-1, 
|G|)=1$.   
Since
$D$ is a Hall subgroup of $G$, it has a complement $M$ in $G$. 
Hence 
$G/E=(D/E)\times (ME/E)$, where
$ME/E\simeq M\simeq G/D$ is $\sigma$-nilpotent. Therefore $G/E$ is
$\sigma$-nilpotent by Lemma 2.2. But then $D\leq E$, a contradiction. Hence $p$
does not divide $|D|$. In particular, $|D|$  is odd.

(8)  {\sl  $D$ is abelian.}

In view of Claim 
(6), $D$ is a Dedekind group.  Hence $D$ is abelian since $|D|$ is  odd  by Claim (7).  

 From Claims (4)--(8) we get that the necessity condition of the theorem holds for $G$.

{\sl Sufficiency.} This directly follows from Theorem B.

The theorem is proved.

\end{document}